\documentclass{amsart} 
\usepackage{amsfonts, amsmath, amssymb, amsthm, color}

\usepackage{graphicx}

\newtheorem{theorem}{Theorem}[section]

\newtheorem{proposition}[theorem]{Proposition}
\newtheorem{lemma}[theorem]{Lemma}

\def\bbm[#1]{\mbox{\boldmath $#1$}}

\newenvironment{Proof}{\noindent{\bf Proof}}{\hfill$\Box$\\[2mm]}

\usepackage{latexsym}
\newcommand{\e}{\epsilon}
\newcommand{\m}{\lambda}

 \renewcommand{\(}{\left(}
\renewcommand{\)}{\right)}
\renewcommand{\[}{\left[}
\renewcommand{\]}{\right]}
\newcommand{\eps}{\epsilon}

\newcommand{\rr}{ \mathbb{R}}



\begin{document}

\title[Steady states with unbounded mass of the Keller-Segel system]{Steady states with unbounded mass of the Keller-Segel system}

\author{Angela Pistoia}
\address[Angela Pistoia] {Dipartimento SBAI, Universt\`{a} di Roma ``La Sapienza", via Antonio Scarpa 16, 00161 Roma, Italy}
\email{pistoia@dmmm.uniroma1.it}
\author{Giusi Vaira}
\address[Giusi Vaira] {Dipartimento SBAI, Universt\`{a} di Roma ``La Sapienza", via Antonio Scarpa 16, 00161 Roma, Italy}
\email{giusi.vaira@sbai.uniroma1.it}

\begin{abstract}
We consider the boundary value problem
\begin{equation*}
\left\{
\begin{array}{lr}
-\Delta u + u =\lambda e^u, \quad \mbox{in}\,\, B_{r_0}\\
\partial_{\nu}u =0\qquad\qquad \mbox{on}\,\, \partial B_{r_0}
\end{array}
\right.
\end{equation*}
where $B_{r_0}$ is the ball of radius $r_0$ in $\mathbb R^N$, $N\geq 2$,
$\lambda>0$ and $\nu$ is the outer normal derivative at $\partial B_{r_0}$. This problem
is equivalent to the stationary Keller-Segel
 system from chemotaxis.\\
We show the existence of a solution  concentrating at the boundary of the ball as $\lambda$ goes to zero.
\end{abstract}


\date{\today}

\maketitle

\section{Introduction}

We consider a system of partial differential equations modelling chemotaxis.
Chemotaxis is a phenomenon of the direct movement of cells in response to the gradient of a chemical, which explains
the aggregation of cells which move towards high concentration of a chemical secreted by themselves.
 The basic model was introduced by Keller and Segel
in \cite{ks} and a simplified form of it reads
\begin{equation}\label{ks-sis}
 \left\{
\begin{array}{lr}
v_t=\Delta v-\nabla\(v\nabla u\) & \hbox{in}\ \Omega\\
\tau u_t=\Delta u-u+v  & \hbox{in}\ \Omega\\
\partial_{\nu}u =\partial_\nu v=0 & \hbox{on}\ \partial\Omega\\
u(x,0)=u_0(x),\ v(x,0)=v_0(x),
\end{array}
\right.
\end{equation}
where $u=u(x,t)\ge0$ and  $v=v(x,t)\ge0$ are the concentration of the species and that of chemical. Here $\Omega$ is a bounded smooth domain in $\rr^N$
and $N\ge2.$ The cases $N=2$ or $N=3$ are of particular interest. In \eqref{ks-sis} $\nu$ denotes the unit outward vector normal at $\partial\Omega$ and $\tau$ is a positive constant.


After the seminal works by Nanjudiah \cite{na} and Childress and Percus \cite{cp} many contributions have been made to the understanding of different analytical aspects of this system and its variations. We refer the reader for instance to \cite{bc,dn,gz,hv1,hv2,hv3,jl,ks,n2,nsy,nss,na,ss1,ss2,ss3,v1,v2,v3}.
\\

In this paper, we study steady states of  \eqref{ks-sis}, namely solutions to the   system
\begin{equation}\label{sis}
 \left\{
\begin{array}{lr}
 \Delta v-\nabla\(v\nabla u\)=0 & \hbox{in}\ \Omega\\
 \Delta u-u+v=0  & \hbox{in}\ \Omega\\
\partial_{\nu}u =\partial_\nu v=0 & \hbox{on}\ \partial\Omega.
\end{array}
\right.
\end{equation}
As point out in \cite{nsy}, stationary solutions to the   Keller-Segel system are of basic importance for the understanding of the global dynamic of the system.

This problem was first studied by Schaaf in \cite{s} in the one dimensional case. In the higher dimensional case Biler in \cite{b} proved the existence of nontrivial radially symmetric solution to \eqref{sis} when $\Omega$ is a ball.
In the general two dimensional case, Wang and Wei in \cite{ww} and Senba and Suzuki in \cite{ss1} proved that    for any $\mu\in \(0,{1\over|\Omega|}+\mu_1\)\setminus \{4\pi m\ :\ m\ge1\}$ problem \eqref{sis} has a non constant solution such that $\int\limits_\Omega v(x)dx=\mu|\Omega|.$
Here  $\mu_1$ is the first eigenvalue of $-\Delta$ with Neumann boundary conditions.
Del Pino and Wei in \cite{dw} reduced system \eqref{sis} to a scalar equation.
Indeed, it is easy to check that $(u,v)$ solves system \eqref{sis} if and only if $v=\lambda e^u$ for some positive constant $\lambda$ and $u$ solves the equation
\begin{equation}\label{p}
 \left\{
\begin{array}{lr}
-\Delta u + u =\lambda e^u & \hbox{in} \ \Omega\\
\partial_{\nu}u =0& \hbox{on} \ \partial\Omega.\\
\end{array}
\right.
\end{equation}
Using this point of view, they proved that for any integers $k$ and $\ell$ there exists a family of solutions $(u_\lambda,v_\lambda)$ to the system \eqref{sis} such that
  $v_\lambda$   exhibits $k$ Dirac measures  inside the domain and $\ell$  Dirac measures on the boundary of the domain as $\lambda\to0,$ i.e.
$$v_\lambda\rightharpoonup\sum\limits_{i=1}^k8\pi\delta_{\xi_i}+\sum\limits_{i=1}^\ell4\pi\delta_{\eta_{i}}\quad\hbox{as}\ \lambda\to0,$$
where $\xi_1,\dots,\xi_k\in\Omega$ and $\eta_1,\dots,\eta_\ell\in\partial\Omega.$ In particular, the solution   has bounded mass, i.e.
$$\lim\limits_{\lambda\to0}\int\limits_{\Omega} v_\lambda(x)dx=\lim\limits_{\lambda\to0}\int\limits_{\Omega} \lambda e^{u_\lambda(x)}dx=4\pi(2k+\ell).$$
 In particular, their argument allows to find a radial solution to the system \eqref{sis}  when $ \Omega$ is a ball in $\rr^2,$ which exhibits a Dirac measure  at the center of the
 ball with mass $8\pi $ when $\lambda$ goes to zero.
  \\

 In the present paper, we find a new radial solution  to the system  \eqref{sis}  when $\Omega$ is a ball in $\rr^N,$ $N\ge2,$ with unbounded mass.
Our main result reads as follows.
\begin{theorem}\label{principale}
Let $\Omega=B(0,r_0) \subset\rr^N,$ $N\ge2,$ be the ball centered at the origin with radius $r_0.$
There exists $\lambda_0$ such that for any $\lambda\in(0,\lambda_0),$ the problem \eqref{p} has a radial solution $(u_{\lambda},v_\lambda)$ such that as $\m\rightarrow 0$
\begin{equation}\label{convergenzaL1}
 \lim\limits_{\lambda\to0}\int\limits_{\Omega} v_\lambda(x)dx=\lim\limits_{\lambda\to0}\int\limits_{\Omega} \lambda e^{u_\lambda(x)}dx=+\infty.
\end{equation}
Moreover, for a suitable choice of positive numbers $\e_\m$ (see \eqref{sceltae}) with $\e_\m\rightarrow 0$ as $\m\rightarrow 0$, we have
\begin{equation}\label{convergenza}
 \lim\limits_{\lambda\to0} \e_\m u_{\m}=\frac{\sqrt{2}}{\mathcal{U}'(r_0)}\mathcal{U} \qquad \hbox{$C^0-$uniformly on compact sets of $\Omega.$}
\end{equation}

\end{theorem}

 Here  $\mathcal{U} $ is the positive radial solution to the problem (see also Lemma \ref{propu})
\begin{equation}\label{G}
\left\{
\begin{array}{lr}
-\Delta \mathcal{U}+\mathcal{U}=0 & \hbox{in} \ B(0,r_0) \\
\mathcal{U} =1& \hbox{on} \ \partial B(0,r_0) .\\
\end{array}
\right.
\end{equation}

 To find this solution, we use a fixed point argument. More precisely, we look for a solution to equation \eqref{p} as $u_\lambda=\bar u_\m+\phi_\m$,
  where the leading term $\bar u_\m$ has to be accurately defined. Once one has a good approximating solution $\bar u_\m$, a simple contraction mapping argument leads to find the higher order term $\phi_\lambda$.

  The difficulty in the construction of the approximated solution $\bar u_\m$ is due to the fact that
$\bar u_\m$ shares the behavior of $\mathcal{U} $   (which solves \eqref{G}) in the inner part of the ball and the behavior of the function $w_\eps$  (see \eqref{w}) near the boundary of the ball.
Here
   \begin{equation}\label{w}
 w_{\e }(r)=\ln\frac{4}{\e^2}\frac{e^{\sqrt{2}\frac{r-r_0}{\e}}}
 {(1+e^{\sqrt{2}\frac{r-r_0}{\e}})^2},\ r\in\rr,\   \e>0.
 \end{equation}
  solve the one dimensional limit problem
 \begin{equation}\label{pblimite}
-w''=e^w\qquad \mbox{in}\,\, \mathbb R .
\end{equation}
  In particular, we have to spend a lot of effort to glue   the   two functions up to the third order (see \eqref{sceltae}, \eqref{v3} and \eqref{closebordo})
  in a neighborhood of the boundary (see Lemma \ref{glue}).\\

  It is important to remark about the analogy existing between our result and some recent results obtained by Grossi in \cite{g1,g2,g3,g4}.
In particular, Grossi and Gladiali in \cite{g4} studied the asymptotic behavior as $\lambda$ goes to zero of the radial solution $z_\lambda$ to the Dirichlet problem
\begin{equation*}
 \left\{
\begin{array}{lr}
-\Delta z   =\lambda e^z & \hbox{in} \ \Omega\\
 z =0& \hbox{on} \ \partial\Omega.\\
\end{array}
\right.
\end{equation*}
when $\Omega=\{x\in\rr^n\ :\ a<|x|<b\}$ is the annulus in $\rr^n.$
In particular, they proved that for a suitable choice of positive numbers $\delta_\m$   with $\delta_\m\rightarrow 0$ as $\m\rightarrow 0,$ $z_\lambda$ satisfies
\begin{equation*}
 \lim\limits_{\lambda\to0} \delta_\m z_{\m}(r)= {2\sqrt{2}} G(r,r^*) \qquad \hbox{$C^0-$uniformly on compact sets of $(a,b).$}
\end{equation*}
where $G(\cdot,r^*)$ is the Green's function of the radial Laplacian with Dirichlet boundary condition and $r^*$ is   suitable choose in $(a,b).$
Moreover, a suitable scaling of $z_\lambda$ in a neighborhood of $r^*$ converges (as $\lambda$ goes to zero) at a solution of the one dimensional limit problem
\eqref{pblimite}.
\\

  The paper is organized as follows.
The definition  of $\bar u_\m$ is given in Section \ref{sec1}, while the construction of a good approximation near the boundary of the ball is carried out in Section \ref{sec2}. In Section \ref{sec3} we estimate the error term and in Section \ref{sec4}  we apply the contraction mapping argument.

\section{The approximated solution}\label{sec1}
We look for a radial solution to the problem \eqref{p}, so we are leading to consider the ODE problem
 \begin{equation}\label{P}
 \left\{
 \begin{array}{lr}
 -u''-\displaystyle\frac{N-1}{r}u'+u=\lambda e^u\qquad \mbox{in}\,\, (0, r_0)\\
 u'(r_0)=0\\
 u'(0)=0.
 \end{array}
 \right.
 \end{equation}
We will construct a solution to \eqref{P} as
$\bar u_{\lambda}+\phi_{\lambda}$ where the leading term $\bar u_\lambda$ is defined as
 \begin{equation}\label{approx}
 \bar u_{\lambda}(r):=\left\{
 \begin{array}{lr}
 u_1(r)\qquad \mbox{in}\,\, (r_0-\delta, r_0)\\
 u_2(r)\qquad \mbox{in}\,\, [r_0-2\delta, r_0-\delta]\\
 u_3(r)\qquad \mbox{in}\,\, (0, r_0-2\delta)
 \end{array}
 \right.
 \end{equation}
and $u_1,$ $u_2$ and $u_3$ are defined as follows.
\\

Basic cells in the construction of the approximate solution $u_1$ near $r_0$ are the functions $w_\eps$     defined in \eqref{w}.
The rate of the concentration parameter $\eps:=\eps_\lambda $ with respect to $\lambda$ is deduced by the relation
\begin{equation}\label{sceltae}
\lambda={4\over\eps_\m^2}e^{-\({{a_1\over\eps_\m}+a_2+a_3\eps_\m}\)},\quad \hbox{i.e.} \quad  \ln\frac{4}{\e_\m^2}-\ln\m=\frac{a_1}{\e_\m}+a_2+a_3\eps_\m
\end{equation}
where  $a_1$, $ a_2$ and $a_3$ are positive constants given in \eqref{abc}. \\
 The right expression  of $u_1 $ is given in \eqref{closebordo}. The construction of $u_1$ is quite involved and it will be carried out in Section \ref{sec2}.

 The approximate solution $u_3$ far away from $r_0$ is build from the function $\mathcal{U} $ which solves \eqref{G} and whose properties are stated in Lemma \ref{propu}.

More precisely,
    \begin{equation}\label{v3}
 u_3(r)=\left(\frac{A_1}{\e_\m}+A_2+A_3\e_\m\right)\mathcal{U}(r)
 \end{equation}
 where  $A_1,$ $ A_2 $ and $A_3$ are positive constants given in \eqref{abc}.\\

 Finally, the approximate solution $u_2$ in the interspace is simply given by
     \begin{equation}\label{u2}
     u_2(r):=\chi(r)u_1(r)+(1-\chi(r))u_3(r) \end{equation}
     where $\chi\in C^2([0,r_0])$ is a cut-off  such that
\begin{equation}\label{chi}\chi\equiv 1\ \hbox{in}\ (r_0-\delta, r_0),\ \chi\equiv 0\ \hbox{in}\ (0, r_0-2\delta),\
|\chi(r)|\leq 1,\   |\chi'(r)|\leq\frac{c}{\delta},\ |\chi''(r)|\leq \frac{c}{\delta^2}.\end{equation}
where the size of the interface $\delta:=\delta_\lambda$ is going to zero with respect to $ \eps$ (or equivalently with respect to $\lambda$)  as
\begin{equation}\label{delta}
\delta_\lambda=\epsilon_\lambda^{\eta},\quad \eta\in \left(\frac 23, 1\right).\end{equation}
The choice of $\eta$ will be made  so that Lemma \eqref{lemerr} holds.
\\

It is important to point out that
 $u_2$   is a good approximation of the solution
 in the interspace, if $u_1$ and $u_3$   perfectly glue in a left neighborhood of $r_0.$
 That implies that we need to go into a third order expansion in $u_1$ (see \eqref{closebordo}) and in $u_3$ (see \eqref{v3}) and also motivates  the rate of $\eps_\m$ made in \eqref{sceltae} and the choice of the constants $A_1,$ $ A_2,$ $ A_3 $ and $a_1$, $a_2,$ $a_3$  made in  Lemma  \ref{glue}.\\\\

  \begin{lemma}\label{propu}
There exists a unique solution to the problem
\begin{equation}\label{grad}
 \left\{
 \begin{array}{lr}
 -\mathcal{U}''-\displaystyle\frac{N-1}{r}\mathcal{U}'+\mathcal{U}=0\qquad \mbox{in}\,\, (0, r_0)\\
 \\
  \mathcal{U}'(0)=0,\ \mathcal{U} (r_0)=1.\\
 \end{array}
 \right.
 \end{equation}
 Moreover
$$0\le \mathcal{U}(r)\le 1\ \hbox{and}\ \mathcal{U}'(r)>0\ \hbox{for any}\ r\in(0,r_0].$$
 \end{lemma}

 \begin{Proof}
The existence and uniqueness of the solution are standard. By the maximum principle we deduce that $\mathcal{U} \le 1$ in $(0,r_0]$.

If $r^*\in (0,r_0)$ is a minimum point of $\mathcal{U}$ with  $\mathcal{U}(r^*)<0,$ by \eqref{grad} we deduce that  $\mathcal{U}''(r^*)=\mathcal{U}(r^*)<0$
which is not possible. So $\mathcal{U} \ge 0$ in $(0,r_0]$.

Finally, we integrate  \eqref{grad} and we get
$$r^{N-1}\mathcal{U} '(r)=\int\limits_0^r t^{N-1}\mathcal{U}  (t)dt\ge0\ \hbox{for any}\ r\in(0,r_0],$$
which implies $\mathcal{U}' >0$ in $(0,r_0]$.
\end{Proof}

 \section{The approximation near the boundary}\label{sec2}
 The function $w_{\epsilon}-\ln\m$ is not a good approximation for our solution near $r_0.$
We will build   some additional correction terms which improve the approximation near $r_0.$
More precisely, we define  the approximation near the point $r_0.$ We define
\begin{equation}\label{closebordo}
u_1(r)=\underbrace{w_{\e }(r)-\ln\m+ \alpha_{\e }(r)}_{1^{st}\ \hbox{order approx.}}+\underbrace{v_{\e }(r)+\beta_{\e }(r)}_{2^{nd}\ \hbox{order approx.}}+\underbrace{z_{\e }(r)}_{3^{rd}\ \hbox{order approx.}}
\end{equation}
where $\alpha_\eps$ is defined in Lemma \ref{lemmaalpha}, $v_\eps$ is defined in Lemma \ref{lemmav},  $\beta_\eps$ is defined in Lemma \ref{lemmabeta}
and $z_\eps$ is defined in Lemma \ref{lemmaz}.
\\\\

 The first term we have to add  is   a sort of projection of the function  $w_{\epsilon},$ namely  the function $\alpha_{\e }$ given in the next  lemma.

  \begin{lemma}\label{lemmaalpha}
\begin{itemize}
\item[(i)] The Cauchy problem
  \begin{equation}\label{alpha}
\left\{
\begin{array}{lr}
-\alpha_{\epsilon, N}''-\displaystyle{\frac{N-1}{r}}\alpha_{\epsilon, N}'=
\displaystyle{\frac{N-1}{r}}w'_{\epsilon}(r)-w_{\epsilon}(r)+\ln\lambda\qquad \mbox{\rm in}\,\, (0, r_0)\\
\\
\alpha_{\epsilon}(r_0)=\alpha_{\epsilon}'(r_0)=0.
\end{array}
\right.
\end{equation}
has the solution
$$\alpha_{\e }(r):=-\int\limits_{r_0}^r{1\over t^{N-1}}\int\limits_{r_0}^t\tau^{N-1}\[{N-1\over \tau}w'_\eps(\tau)-w_\eps(\tau)+\ln\lambda\]d\tau\ dt.$$
\item[(ii)]  The following expansion holds
\begin{equation}\label{alfa}
\alpha_{\e }(\e s+r_0)= \e\alpha_1(s)+\e^ 2\alpha_2(s)+
 O\left(\e^3s^4\right)
\end{equation}
where
\begin{equation}\label{alfa1}
\alpha_{1 }(s):= - \frac{N-1}{r_0}\int_0^s w(\sigma)\ d\sigma  +\frac12 a_1 s^2
\end{equation}
and
\begin{eqnarray}\label{alfa2}
\alpha_{2 }(s)&:= &  \int_0^s\int_0^\sigma \[w (\rho)-\ln4\]\, d\rho\ d\sigma +{(N-1)(N-2)\over r_0^2}\int_0^s\int_0^\sigma w(\rho)\, d\rho\ d\sigma\nonumber\\ & &+{(N-1)\over r_0 ^2 }\int_0^s\sigma w(\sigma)\ d\sigma
-{N-1\over6r_0}  a_1 s^3+\frac12 a_2 s^2
\end{eqnarray}
\item[(iii)]
For any $r\in (0, r_0-\delta)$
\begin{eqnarray}\label{alfa3}
\alpha_{\e }(r)&=&-\frac{(N-1)\ln 4}{r_0}(r-r_0)+\left[\frac{(N-1)^2\ln 4}{r_0^2}
-\frac{\sqrt{2}(N-1)}{\e r_0}+\ln\frac{4}{\e^2}-\ln\m\right]\frac{(r-r_0)^2}{2}\nonumber\\
&&+\left[\frac{N(N-1)\sqrt{2}}{\e r_0^2}+\frac{\sqrt{2}}{\e}-\frac{N-1}{r_0}\left(\ln\frac{4}{\e^2}
-\ln\m\right)\right]\frac{(r-r_0)^3}{6}\nonumber\\
&&+O\left(\frac{(r-r_0)^4}{\e}\right)+O\left((r-r_0)^3\right)
\end{eqnarray}
 \end{itemize}

 \end{lemma}
 \begin{Proof}

{\em Proof of (i).}
It is just a straightforward computation.\\

{\em Proof of (ii).}

We get (setting $t=\eps\sigma+r_0$ and $\tau=\eps\rho+r_0$)
\begin{eqnarray*}
\alpha_{\e }(\e s+r_0)& &= -\eps^2 \int\limits_{0}^s{1\over (\eps \sigma+r_0)^{N-1}}\int\limits_{0}^ \sigma(\eps\rho+r_0)^{N-1}\[{N-1\over\eps\rho+r_0}{1\over\eps}w '(\rho)-\[w(\rho)-\ln4\]+\ln\lambda-\ln{4\over\eps^2}\]d\sigma\ d\rho\\
& &=-\eps^2 \int\limits_{0}^s\({1\over r_0^{N-1}}-{{N-1\over r_0^N}\eps\sigma}\)
 \int\limits_{0}^ \sigma\(r_0^{N-1}+(N-1)r_0^{N-2}\eps\rho\)\times
\\ & &\qquad\times
\[(N-1)\({1\over r_0}-{{ 1\over r_0^2}\eps\rho}\){1\over\eps}w '(\rho)-\[w(\rho)-\ln4\]+\ln\lambda-\ln{4\over\eps^2}\]d\sigma\ d\rho\\
& &+O\left(\e^3s^4\right)
\end{eqnarray*}
Here we used that
$$w_\eps(r)=\ln{4\over\eps^2}+ w\({r-r_0\over\eps}\)-\ln4 \ \hbox{and}\ w'_\eps(r)={1\over\eps}w'\({r-r_0\over\eps}\).$$
The claim follows by \eqref{sceltae}.
\\

{\em Proof of (iii).}
Set $\bar w_\eps(r):=w_\eps(r)-\ln{1\over\eps^2}.$

We have
\begin{eqnarray*}
\alpha_{\e }(r)& &=-\int\limits_{r_0}^r{1\over t^{N-1}}\int\limits_{r_0}^t\tau^{N-1}\[{N-1\over \tau}w'_\eps(\tau)-w_\eps(\tau)+\ln\lambda\]d\tau\ dt\\
& &= -\int\limits_{r_0}^r{1\over t^{N-1}}\int\limits_{r_0}^t\tau^{N-1}\[{N-1\over \tau}\bar w'_\eps(\tau)-\bar w_\eps(\tau)+\(\ln\lambda-\ln{1\over\eps^2}\)\]d\tau\ dt\\
&&=-(N-1)\int\limits_{r_0}^r{\bar w_\eps(t)\over t}dt+\int\limits_{r_0}^r{1\over t^{N-1}}\int\limits_{r_0}^t\[(N-1)(N-2)\tau^{N-3}+\tau^{N-1}\]\bar w_\eps(\tau)d\tau\ dt\\
& &+\(\ln\lambda-\ln{1\over\eps^2}\)
\left\{
\begin{array}{lr}{1\over 2N}(r_0^2-r^2)+{r_0^2\over 2}\ln{r\over r_0}\ \hbox{if}\ N=2\\
\\
 {1\over 2N}(r_0^2-r^2)+{r_0^N\over N(N-2)}\({1\over r_0^{N-2}}-{1\over r^{N-2}}\)\ \hbox{if}\ N\ge3.\\
\end{array}
\right.
\end{eqnarray*}

Now we observe that
in $[r_0-2\delta, r_0-\delta]$ we have
\begin{equation}\label{log}
\ln\frac{r}{r_0}=\ln\left(1+\frac{r-r_0}{r_0}\right)=\frac{r-r_0}{r_0}-\frac{(r-r_0)^2}{2r_0^2}+\frac{(r-r_0)^3}{3r_0^3}+O\left((r-r_0)^4\right)
\end{equation}
\begin{equation}\label{rN}
\frac{1}{r^{N-2}}=\frac{1}{r_0^{N-2}}-\frac{(N-2)}{r_0^{N-1}}(r-r_0)+\frac{(N-2)(N-1)}{r_0^N}\frac{(r-r_0)^2}{2}
-\frac{N(N-1)(N-2)}{r_0^{N+1}}\frac{(r-r_0)^3}{6}+O((r-r_0)^4)
\end{equation}
and also
\begin{equation}\label{expbar}
\bar{w}_{\e }(s)=\ln 4 +\frac{\sqrt{2}}{\e}(s-r_0)+O\left(e^{-\frac{|s-r_0|}{\e}}\right).
 \end{equation}

A tedious but straightforward computation proves our claim.

\end{Proof}
\\\\

The function $w_{\e }(r)-\ln\m +\alpha_{\e }(r)$ is yet   a bad approximation of the solution near the boundary point $r_0$.
We have to add a correction term $v_\eps$ given in next lemma, which solves a linear problem and {\em kills} the $\eps-$order term in \eqref{alfa}.

\begin{lemma}\label{lemmav}
\begin{itemize}
\item[(i)]
There exists a solution $v$ of the linear problem (see \eqref{alfa1})
\begin{equation}\label{VV}
-v''-e^w v= e^w\alpha_1\quad\hbox{in}\ \rr
\end{equation}
such that
$$v(s)=\nu_1s+\nu_2+O\(e^s\)\quad{and}\quad v'(s)=\nu_1 +O\(e^s\)\quad\hbox{as}\ s\to-\infty $$
where $\nu_2\in\rr$ and
\begin{equation}\label{nu1}\nu_1:= -\frac{2(N-1)}{r_0}(1-\ln 2) +a_1{\sqrt{2}}  { \ln 2  }  \end{equation}
\item[(ii)] In particular, the function  $v_{\e }(r):=\e v\left(\frac{r-r_0}{\e}\right) $ is a solution of the linear problem
\begin{equation}\label{ve}
-v_{\e }''-e^{w_{\e}}v_{\e }=\eps e^{w_{\e}(r)} {\alpha_1} \(\frac{r-r_0 }{\eps}\right)\quad\hbox{in}\ \rr
\end{equation}
such that if $r\in[0,r_0-\delta]$ it satisfies
\begin{equation}
\label{vepsi}
v_{\e }(r)=\nu_1(r-r_0)+\nu_2 \e +O(\e e^{-\frac{|r-r_0|}{\e}})\quad{and}\quad v'_{\e }(r)=\nu_1 O(e^{-\frac{|r-r_0|}{\e}})\quad\hbox{as}\ \e\rightarrow 0.\end{equation}

\end{itemize}
\end{lemma}
\begin{Proof}
The  result immediately follows by Lemma \ref{YY}.
In our case is
    $$\nu_1:={1\over\sqrt2}\int\limits_{-\infty}^0 \(-\frac{N-1}{r_0}\int_0^r w(y)\, dy+a_1\frac{r^2}{2}\)w'(r)e^{w}(r)dr$$
   and a straightforward computation proves \eqref{nu1}.
  \end{Proof}

  \begin{lemma}\label{YY} [\cite{g3}, Lemma 4.1]
 Let $h:\rr\to\rr$ be a continuous function. The function
  \begin{equation}\label{esprY}
  Y(t)=w'(t)\int_0^t\frac{1}{w'(s)^2}\left(\int_s^0 h(z)w'(z)e^w\, dz\right)\, ds
  \end{equation}
  is a solution to
  \begin{equation}\label{Y}
 -Y''-e^{w}Y=he^w\qquad \mbox{in}\,\, \rr
 \end{equation}
Moreover, it satifies
\begin{eqnarray*}
&&Y(t)=\frac{t}{\sqrt{2}}\int_{-\infty}^0h(r)w'(r)e^{w}\,dr-\int_{-\infty}^0 \left(\frac{2}{1-e^{\sqrt{2}s}}+\frac{s}
{\sqrt{2}}\right)h(s)w'(s)e^w\, ds+O\left(e^t\right),\\
&&Y'(t)=\frac{1}{\sqrt{2}}\int_{-\infty}^0h(r)w'(r)e^{w}\,dr+O\left(e^t\right)\quad \hbox{as}\ t\rightarrow-\infty
\end{eqnarray*}
  and
\begin{eqnarray*}
&&Y(t)=\frac{t}{\sqrt{2}}\int_0^{+\infty}h(r)w'(r)e^{w}\,dr-\int_0^{+\infty} \left(\frac{2}{1-e^{\sqrt{2}s}}+\frac{s}
{\sqrt{2}}\right)h(s)w'(s)e^w\, ds+O\left(e^{-t}\right)\\
& &Y'(t)=\frac{1}{\sqrt{2}}\int_0^{+\infty}h(r)w'(r)e^{w
}\,dr+O\left(e^{-t}\right)\quad \hbox{as}\ t\rightarrow+\infty.
\end{eqnarray*}
 \end{lemma}

As we have done for the function $w_\eps,$ we have to add  the  projection   of the function  $v_{\epsilon },$ namely  the function $\beta_{\e }$ given in the next  lemma.

 \begin{lemma}\label{lemmabeta}
 \begin{itemize}\item[(i)]
The  Cauchy problem:
\begin{equation}\label{beta}
\left\{
\begin{array}{lr}
-\beta_{\e }''-\displaystyle\frac{(N-1)}{r}\beta_{\e }' = \frac{(N-1) }{r}v'_\eps(r)\\\\
\beta_{\e }(r_0)=\beta_{\e }'(r_0)=0.
\end{array}
\right.
\end{equation}
has the solution
$$\beta_{\e }(r)=-(N-1)\int_r^{r_0}\frac{1}{t^{N-1}}\int_t^{r_0}\tau^{N-2}v'_{\e }(\tau)\, d\tau\, d\, t.$$
\item[(ii)]  The following expansion holds
\begin{equation}\label{beta2}
\beta_{\e }(\e s+r_0)=
  \e^2\beta_1(s)+O\(\eps^3s^3\),\quad \beta_1(s):=-\frac{(N-1)}{r_0}\int_0^s\int_0^\sigma v'(\rho)\, d\rho\,d\sigma.\end{equation}
\item[(iii)]
For any $r\in (0, r_0-\delta)$
\begin{equation}\label{beta1}
\beta_{\e }(r)=-\frac{(N-1)\nu_1}{r_0}\frac{(r-r_0)^2}{2}+O((r-r_0)^3)\end{equation}
\end{itemize}
\end{lemma}
\begin{Proof}
We argue as in Lemma \ref{lemmaalpha}.
\end{Proof}

Unfortunately, the function $w_{\e, r_0}(r)-\ln\m +\alpha_{\e }(r)+v_\eps(r)+\beta_\eps(r)$ is yet   a bad approximation of the solution near the boundary point $r_0$.
We have to add an extra correction term $z_\eps$ given in next lemma, which solves a linear problem and {\em kills} all the $\eps^2-$order terms (in particular, those in \eqref{alfa} and in \eqref{beta2}).

\begin{lemma}\label{lemmaz}

\begin{itemize}
\item[(i)]
There exists a solution $z$ of the linear problem  (see \eqref{alfa1}, \eqref{alfa2}, \eqref{beta2}, \eqref{VV})
\begin{equation}\label{ZZ}
-z''-e^w z= e^z\[\alpha_2(s)+\beta_1(s)+{1\over2}\(\alpha_1(s)+v(s)\)^2\]\quad\hbox{in}\ \rr
\end{equation}
such that
$$z(s)=\zeta_1s+\zeta_2+O\(e^s\)\quad{and}\quad z'(s)=\zeta_1 +O\(e^s\)\quad\hbox{as}\ s\to-\infty $$
where $\zeta_1,\zeta_2\in\rr.$
\item[(ii)] In particular, the function  $z_{\e }(r):=\e^2 z\left(\frac{r-r_0}{\e}\right) $ is a solution of the linear problem
\begin{equation}\label{z}
-z_{\e }''-e^{w_{\e }}z_{\e }= \eps^2 e^{w_{\e }}
\left\{\alpha_2\({r-r_0\over\eps}\)+ \beta_1\({r-r_0\over\eps}\)
 +{1\over2}\[ \alpha_1\({r-r_0\over\eps}\)+v\({r-r_0\over\eps}\)\]^2\right\}
\end{equation}
such that if $r\in[0,r_0-\delta]$ it satisfies
\begin{equation}
\label{zepsi}
z_{\e }(r)=\eps\zeta_1(r-r_0)+\zeta_2 \e^2 +O\(\e^2 e^{-\frac{|r-r_0|}{\e}} \)\quad\hbox{as}\ \e\rightarrow 0.\end{equation}

\end{itemize}

\end{lemma}
\begin{Proof}

 The  result immediately follows by Lemma \ref{YY}.

\end{Proof}

\section{The error estimate  }\label{sec3}

Let us define the error term
\begin{equation}\label{rlam}\mathcal{R}_{\m}(\bar u_\lambda)=-\bar u_\lambda''-\frac{N-1}{r}\bar u'_\lambda+\bar u_\lambda-\m e^{\bar u_\lambda}.\end{equation}
where $\bar u_\lambda$ is defined as in \eqref{approx}.

First of all, it is necessary to  choose constants $a,$ $b$ and $c$ in \eqref{sceltae} and $A_1,$ $A_2$ and $A_3$ in \eqref{v3}
such that the approximate solutions in
the neighborhood of the boundary   and inside the interval glue up.

\begin{lemma}\label{glue}
If
\begin{equation}\label{abc}
a_1=A_1:={\sqrt 2\over \mathcal{U}'(r_0)},\ a_2 =A_2:={1\over \mathcal{U}'(r_0)}\({\ln 4\over \mathcal{U}'(r_0)}-2{N-1\over r_0}\),\ a_3 :=A_3-\nu_2,\ A_3:={\zeta\over\mathcal{U}'(r_0)}
\end{equation}
then for any $r\in[r_0-2\delta,r_0-\delta]$ we have
\begin{eqnarray*}
& & u_1(r)-u_3(r)=O\(e^{-{|r-r_0|\over\eps}}\)+O\(\eps^2\)+ O\(\eps{(r-r_0)^2 }\)+O\({(r-r_0)^3 }\)+O\({(r-r_0)^4\over\eps}\),\\
& & u'_1(r)-u'_3(r)=O\({1\over\eps}e^{-{|r-r_0|\over\eps}}\)+O\(\eps \)+ O\(\eps{(r-r_0)  }\)+O\({(r-r_0)^2 }\)+O\({(r-r_0)^ 3\over\eps}\).
\end{eqnarray*}
\end{lemma}

\begin{Proof} Let us prove the first estimate. The proof of the second estimate is similar.
By \eqref{sceltae}, by \eqref{alfa3}, \eqref{vepsi}, \eqref{beta1} and \eqref{zepsi}
we deduce that if $r\in[r_0-2\delta,r_0-\delta]$ then
\begin{eqnarray}\label{ucon1}
u_1(r)=& &\[\ln{4\over\eps^2}-\ln\lambda+\nu_2\eps\]+\[{\sqrt2\over\eps}-{(N-1) \ln4\over r_0}+\nu_1+\zeta_1\eps\](r-r_0)\nonumber\\
& &+\[{(N-1)^2 \ln4\over r_0^2}-{\sqrt2(N-1)\over r_0}{1\over\eps}+\ln{4\over\eps^2}-\ln\lambda-{\nu_1(N-1)\over r_0}\]{(r-r_0)^2\over2}\nonumber\\
 & &+\[{N(N-1) \sqrt2\over r_0^2}{1\over\eps}+{\sqrt2(N-1)}{1\over\eps}-{N-1\over r_0}\(\ln{4\over\eps^2}-\ln\lambda\)\]{(r-r_0)^3\over6}\nonumber\\
& &+O\(e^{-{|r-r_0|\over\eps}}\)+O\(\eps^2\)+O\({(r-r_0)^4\over\eps}\)+O\({(r-r_0)^3 }\)\nonumber\\
& &=\[{a_1\over\eps}+a_2+a_3\eps+\nu_2\eps\]+\[{\sqrt2\over\eps}-{2(N-1)  \over r_0}+a_1\sqrt2\ln2+\zeta_1\eps\](r-r_0)\nonumber\\
& &+\[-{(N-1)\sqrt2\over r_0}{1\over\eps}+{a_1\over\eps}+a_2+2{(N-1)^2 \over r_0^2}-{a_1(N-1)\sqrt2\ln2\over r_0} \]{(r-r_0)^2\over2}\nonumber\\
 & &+\[ {N(N-1) \sqrt2\over r_0^2}{1\over\eps}+{\sqrt2  \over\eps}-{a_1(N-1)\over r_0} {1\over \eps}\]{(r-r_0)^3\over6}\nonumber\\
& &+O\(e^{-{|r-r_0|\over\eps}}\)+O\(\eps^2\)+O\({(r-r_0)^4\over\eps}\)+O\({(r-r_0)^3 }\)
\end{eqnarray}

On the other hand, by the mean value Theorem we deduce that
$$\mathcal{U}(r)=\mathcal{U}(r_0)+\mathcal{U}'(r_0)(r-r_0)+\mathcal{U}''(r_0){(r-r_0)^2\over2}+\mathcal{U}''''(r_0){(r-r_0)^3\over6}+O\({(r-r_0)^4 }\)$$
with $\mathcal{U}(r_0)=1,$
$$ \mathcal{U}''(r_0)=-{N-1\over r_0}\mathcal{U}'(r_0)+\mathcal{U}(r_0)=-{N-1\over r_0}\mathcal{U}'(r_0)+1$$
and
$$\mathcal{U}'''(r_0)=-{N-1\over r_0}\mathcal{U}''(r_0)+ {N-1\over r_0^2}\mathcal{U}'(r_0)+\mathcal{U}'(r_0)={N(N-1)\over r_0^2}\mathcal{U}'(r_0)+\mathcal{U}'(r_0)-{N-1\over r_0}.$$
These relations easily follow by differentiating \eqref{grad}.
Therefore, if $r\in[r_0-2\delta,r_0-\delta]$   we have
\begin{eqnarray}\label{ucon3}
u_3(r)= \({A_1\over\eps}+A_2+A_3\eps\)\mathcal{U}(r)= & &\({A_1\over\eps}+A_2+A_3\eps\)+ \({A_1\over\eps}+A_2+A_3\eps\)\mathcal{U}'(r_0)(r-r_0)\nonumber\\
& &+\mathcal{U}''(r_0)\({A_1\over\eps}+A_2\){(r-r_0)^2\over2}+\mathcal{U}'''(r_0 ){A_1\over\eps} {(r-r_0)^3\over6}\nonumber\\
& &  +O\(\eps(r-r_0)^2\)+O\({(r-r_0)^3 }\)+O\({(r-r_0)^4\over\eps}\)
\end{eqnarray}
If \eqref{abc} holds then combining \eqref{ucon1} and \eqref{ucon3} we easily get the claim.

\end{Proof}

\begin{lemma} \label{lemerr} There exists $C>0$ and $\lambda_0>0$ such that for any $\m\in(0,\m_0)$ we have
$$\|\mathcal{R}_{\m}\|_{L^1(0,r_0) }=O\(\e_\m^{ 1+\sigma} \)\quad \hbox{for some}\ \sigma>0.
$$
 \end{lemma}

\begin{Proof}

{\em Step 1. Evaluation of the error in $(r_0-\delta, r_0)$.
}\\
We use this estimate
 $1-e^t=-t-\frac{t^2}{2}+O(t^3)$ and we get
\begin{eqnarray*}
 \mathcal{R}_{\m}(u_1)&=& -u_1''-\frac{N-1}{r}u_1'+u_1-\m e^{u_1}\\
& =& -w_{\e}''-\frac{N-1}{r_0}w_{\e}'+w_{\e}-\ln\m -\alpha_{\e }''-\frac{N-1}{r}\alpha_{\e }'\\
&& +\alpha_{\e }-v_{\e }''-\frac{N-1}{r}v_{\e }'+v_{\e }-\beta_{\e }''-\frac{N-1}{r}\beta_{\e }'+\beta_{\e }\\
&&  -z_{\e }''-\frac{N-1}{r}z_{\e }'+z_{\e }-\m e^{w_{\e, r_0}-\ln\m +\alpha_{\e }+v_{\e }+\beta_{\e }+z_{\e }} \\
&=& \alpha_{\e}+v_{\e}+\beta_\eps+z_\eps-\frac{N-1}{r}z_{\e}'\\ & &+
  e^{w_{\e }}\left\{1-e^{\alpha_{\e}+v_{\e}+\beta_{\e}+z_{\e}}
+v_{\e}+z_{\e} +  \eps  \alpha_1\({r-r_0\over\eps}\) \right.\\
&&  +  \eps^2\left.\[\alpha_2\({r-r_0\over\eps}\)+ \beta_1\({r-r_0\over\eps}\)
 +{1\over2}\( \alpha_1\({r-r_0\over\eps}\)+v\({r-r_0\over\eps}\)\)^2\]\right\}\\
 &=& \alpha_{\e}+v_{\e}+\beta_\eps+z_\eps-\frac{N-1}{r}z_{\e}'\\ & &+
  e^{w_{\e }}\left\{ -\alpha _\eps  -\beta_\eps-{1\over2} \( \alpha_{\e}+v_{\e} \)^2+\eps\alpha_1\({r-r_0\over\eps}\) \right.\\
&&  +  \eps^2\left.\[\alpha_2\({r-r_0\over\eps}\)+ \beta_1\({r-r_0\over\eps}\)
 +{1\over2}\( \alpha_1\({r-r_0\over\eps}\)+v\({r-r_0\over\eps}\)\)^2\]\right\}\\
 &&+O\(e^{w_{\e }}|{\alpha_{\e}+v_{\e}+\beta_{\e}+z_{\e}}|^3\)+O\(e^{w_{\e }}|{  \beta_{\e}+z_{\e}}|^2\)
 +O\(e^{w_{\e }}|\({  \alpha_{\e}+v_{\e}}\)\( {  \beta_{\e}+z_{\e}} \)|\)\\
\end{eqnarray*}
because $\alpha_\eps$ solves \eqref{alpha}, $v_\eps$ solves \eqref{ve}, $\beta_\eps$ solves \eqref{beta} and $z_\eps$ solves \eqref{z}.

We have
$$\int\limits_{r_0-\delta}^{r_0}|{\alpha_{\e}+v_{\e}+\beta_{\e}+z_{\e}}| (r) dr=O\(\int\limits_{r_0-\delta}^{r_0} {(r-r_0)^2\over\eps} dr\)=O\({\delta^3\over\eps}\)=O\(\eps^{3\eta-1}\),
$$
because by  \eqref{alfa}, \eqref{beta2}, the properties of $v_\eps$ in Lemma \ref{lemmav}  and $z_\eps$   in \ref{lemmaz}
we deduce
$$ \alpha_{\e}(r)=O\({(r-r_0)^2\over\eps}\),\  \beta_{\e}(r)=O\((r-r_0)^2 \),\ v_\eps(r)=O\(|r-r_0|+\eps\),\ z_\eps(r)=O\(\eps|r-r_0|+\eps^2\).$$

By Lemma \ref{lemmaz} we also deduce that $z_{\e}'(r)=O(\eps)$ and so
$$
\int\limits_{r_0-\delta}^{r_0}\left|\frac{1}{r}z_{\e}'(r)\right| dr=O\(\eps\delta\)=O\(\eps^{ 1+\eta}\).
$$

Moreover, we scale $s=\eps r+r_0$ and we get
 \begin{eqnarray*}
 &  &\int\limits_{r_0-\delta}^{r_0}  e^{w_{\e }}\left| -\alpha _\eps  -\beta_\eps-{1\over2} \( \alpha_{\e}+v_{\e} \)^2+\eps\alpha_1\({r-r_0\over\eps}\) \right.\\
&&  +  \eps^2\left.\[\alpha_2\({r-r_0\over\eps}\)+ \beta_1\({r-r_0\over\eps}\)
 +{1\over2}\( \alpha_1\({r-r_0\over\eps}\)+v\({r-r_0\over\eps}\)\)^2\]\right|dr\\
 & =&{1\over\eps}\int\limits_{ -\delta/\eps}^{ 0}  e^{w(s)}\left| -\alpha _\eps (\eps s+r_0) -\beta_\eps (\eps s+r_0)-{1\over2} \( \alpha_{\e} (\eps s+r_0)+\eps v(s) \)^2+\eps\alpha_1\(s\) \right.\\
&&  +  \eps^2\left.\[\alpha_2\(s\)+ \beta_1\(s\)
 +{1\over2}\( \alpha_1\(s\)+v\(s\)\)^2\]\right|ds\\
 &=  &O\(\eps^2 \int\limits_{ \rr}  e^{w(s)}s^3\, ds\)=O\(\eps^2  \).\end{eqnarray*}

 Finally, we scale $s=\eps r+r_0$  and we get
 \begin{eqnarray*}
 & &\int\limits_{r_0-\delta}^{r_0}e^{w_{\e }}|{\alpha_{\e}+v_{\e}+\beta_{\e}+z_{\e}}|^3dr=O\(\int\limits_{r_0-\delta}^{r_0}e^{w_{\e }}\(|{\alpha_{\e}|^3+|v_{\e}|^3+|\beta_{\e}|^3+|z_{\e}}|^3\)dr\)\\ & &=O\(\eps^2 \int\limits_{ \rr}  e^{w(s)}s^6ds
 +\eps^2 \int\limits_{ \rr}  e^{w(s)}v^3(s)ds+\eps^5 \int\limits_{ \rr}  e^{w(s)}s^6ds+\eps^5 \int\limits_{ \rr}  e^{w(s)}z^3(s)ds \)=O\(\eps^2  \)\\
&&\int\limits_{r_0-\delta}^{r_0}e^{w_{\e }}|{ \beta_{\e}+z_{\e}}|^2dr=O\(\eps^3 \int\limits_{ \rr}  e^{w(s)}s^4\, ds+\eps^3 \int\limits_{ \rr}  e^{w(s)}z^2\, ds\)=O\(\eps^3  \)\\
&&\int\limits_{r_0-\delta}^{r_0}e^{w_{\e }}|\({  \alpha_{\e}+v_{\e}}\)\( {  \beta_{\e}+z_{\e}} \)|dr=
O\(\eps^2 \int\limits_{ \rr}  e^{w(s)}(s^2+|v|)(s^2+|z|)\, ds\)=O\(\eps^2  \),
 \end{eqnarray*}
because by  \eqref{alfa} and  \eqref{beta2}
we deduce
$$ \alpha_{\e}(\eps s+r_0)=O\({\eps s^2 }\),\  \beta_{\e}(\eps s+r_0)=O\(\eps^2s^2 \).$$

By collecting all the previous estimates and taking into account the choice of $\eta$ in \eqref{delta} we get
\begin{equation}\label{errore1}
\|\mathcal{R}_{\m}\|_{L^1(r_0-\delta, r_0)}=O\(\e^{ 1+\sigma}\)\quad \hbox{for some}\ \sigma>0.
\end{equation}

\medskip

{\em Step 2: Evaluation of the error in $(0, r_0-2\delta)$}.\\

First of all,    if $\delta$ is small enough (namely $\eps$ is small enough) we have
 $$\mathcal{U}(r)\le \mathcal{U}(r_0-2\delta)=\mathcal{U}(r_0)+\mathcal{U}'(r_0)(-2\delta)+{1\over2}\mathcal{U}''(  r_0-2\theta\delta) (2\delta)^2 \le
 1 -2\mathcal{U}'(r_0) \delta .$$
 because $\mathcal{U}$ is increasing (see Lemma \ref{propu}) and the mean value theorem applies for some $\theta\in(0,1).$

Therefore, by \eqref{sceltae}, \eqref{delta} and \eqref{abc}, we get
\begin{eqnarray*}
& &\mathcal{R}_{\m}(u_3)= -u_3''-\frac{N-1}{r}u_3'+u_3-\m e^{u_3}=
-\m e^{\left(\frac{A_1}{\e}+A_2+A_3\e\right)\mathcal{U}(r)}=-{4\over\eps^2}e^{(A_3-a_3)\eps} e^{\left(\frac{A_1}{\e}+A_2+A_3\e\right)\[\mathcal{U}(r)-1\]}\\
& &=O\({1\over\eps^2}e^{-2A_1\mathcal{U}'(r_0)\frac{\delta}{\e} }\) =O\({1\over\eps^2}e^{-2\sqrt2\frac{1}{\e^{1-\eta}} }\).
 \end{eqnarray*}

This implies that
\begin{equation}\label{errore2}
\|\mathcal{R}_{\m}(u_3)\|_{L^1(0,r_0-2\delta)}=O\(\e^{ 1+\sigma}\)\quad \hbox{for any}\ \sigma>0.
\end{equation}

{\em Step 3: Evaluation of the error in $[r_0-2\delta, r_0-\delta]$}\\
We recall that $u_2=\chi u_1+(1-\chi) u_3$ hence
\begin{eqnarray*}
\mathcal{R}_{\m}(u_2)&=& \chi\left[-u_1''-\frac{N-1}{r}u_1'+u_1\right]+(1-\chi)
\left[-u_3''-\frac{N-1}{r}u_3'+u_3\right]\\
&& -2\chi'\left(u_1'-u_3'\right)+\left[-\chi''-\frac{N-1}{r}\chi'+\chi\right](u_1-u_3)-\lambda e^{\chi(u_1-u_3)+u_3}\\
&=& \chi \mathcal{R}_{\m}(u_1)+(1-\chi)\mathcal{R}_{\m}(u_3)-\m \chi e^{u_1} \[e^{(\chi-1)(u_1-u_3)}-1\] +\lambda(1-\chi) e^{u_3}\\
&&-2\chi'\left(u_1'-u_3'\right)+\left[-\chi''-\frac{N-1}{r}\chi'+\chi\right](u_1-u_3)
\end{eqnarray*}
By Lemma \eqref{glue} we immediately get (taking into account the choice of $\eta$ in \eqref{delta})
$$\int\limits_{r_0-2\delta}^{r_0-\delta}\left|\chi'(r)\left(u_1'(r)-u_3'(r)\right)\right|dr=O\(\delta^2\)=O\(\eps^{1+\sigma}\),$$
$$\int\limits_{r_0-2\delta}^{r_0-\delta}\left|\left[-\chi''(r)-\frac{N-1}{r}\chi'(r)+\chi(r)\right](u_1(r)-u_3(r))\right|(r)dr=O\(\delta^2\)=O\(\eps^{1+\sigma}\),$$
and
$$\int\limits_{r_0-2\delta}^{r_0-\delta}\left|\m \chi e^{u_1(r)} \[e^{(\chi(r)-1)(u_1(r)-u_3(r))}-1\] \right| dr=O\(\int\limits_{r_0-2\delta}^{r_0-\delta}\lambda
  e^{u_1(r)} |u_1(r)-u_3(r)|   dr \)=O\(\lambda\eps^2\),$$
  because $e^t-1=O(t).$
Arguing exactly as in Step 1 one proves that
$$\int\limits_{r_0-2\delta}^{r_0-\delta}\chi(r)\left| \mathcal{R}_{\m}(u_1)(r)\right|dr=O\(\eps^{1+\sigma}\) $$
and arguing exactly as in   Step 2   one proves that
$$\int\limits_{r_0-2\delta}^{r_0-\delta}\(1-\chi(r)\)\left| \mathcal{R}_{\m}(u_3)(r)\right|dr=O\(\eps^{1+\sigma}\)\quad\hbox{and}\
\int\limits_{r_0-2\delta}^{r_0-\delta}\lambda\(1-\chi(r)\) e^{u_3}(r)dr=O\(\eps^{1+\sigma}\).$$

Collecting all the previous estimates, we get
\begin{equation}\label{errore3}
\|\mathcal{R}_{\m}(u_2)\|_{L^1(r_0-2\delta,r_0- \delta)}=O\(\e^{ 1+\sigma}\)\quad \hbox{for some}\ \sigma>0.
\end{equation}

\medskip
The claim follows by \eqref{errore1}, \eqref{errore2} and \eqref{errore3}.

\end{Proof}

\begin{lemma}\label{conv}
It holds that
\begin{equation}\label{conv1}
\m \e_\lambda^2 e^{u_{\m}(\e_\lambda s+r_0)}\rightarrow e^{w(s)}\quad  \hbox{$C^0-$uniformly on compact sets of $(-\infty, 0]$ as}\ \lambda\rightarrow 0
\end{equation}
and
\begin{equation}\label{conv2}
\m \e_\lambda  \int\limits_0^{r_0}e^{u_{\m}(r)}dr\rightarrow    \int\limits_ {\rr}e^{w(s)}ds\quad  \hbox{  as}\ \lambda\rightarrow 0\end{equation}
\end{lemma}
\begin{Proof}
 Let $[a, b]\subset (-\infty, 0]$. If  $\lambda$ is small enough then
$$u_{\m}(\e_\lambda s+r_0)=u_1(\e_\lambda s+r_0)\quad \hbox{for any}\ s\in [a,b].$$

On the other hand, by   \eqref{alfa}, \eqref{beta2}, the properties of $v_\eps$ in Lemma \ref{lemmav}  and $z_\eps$   in \ref{lemmaz}
we deduce
$$ \alpha_{\e}(\e  s+r_0)+\e  v(s)+\beta_{\e}(\e  s+r_0)+\e ^2 z(s)= O\({\eps^2}\)+O\({\eps  |s|+\eps }\) +O\({\eps  ^2 s^2}\)+O\({\eps^2  |s|+\eps^2 }\) $$
and so
 $$u_1  (\e  s+r_0)=w(s)+\ln\frac{1}{\e^2}-\ln \m +O\({\delta  |s|+\delta }\).$$
Therefore,
\begin{equation}\label{cc1}\m \e_\lambda^2 e^{u_{\m}(\e_\lambda s+r_0)}=e^{w(s)+O\({\delta  |s|+\delta }\)}\end{equation}
and \eqref{conv1} follows, since $s\in[a,b].$

   Moreover, since
  $w(s)=\sqrt 2s+O\(e^{\sqrt 2s}\)$ as $s$ goes to $-\infty,$ we also deduce that
  if $\m$ (and also $\delta$) is small enough  there exist   $a,b>0$ such that
\begin{equation}\label{cc2}
\m \e^2 e^{u _1(\e  s+r_0)}\leq be^{-a|s|}\quad\hbox{for any}\ s\in (-\infty,0].\end{equation}

Now, we have (scaling $r=\eps s+r_0$ in the first integral and arguing as in Step 3 of Lemma \ref{lemerr} to estimate the second and the third integral)
\begin{eqnarray*}\m \e_\lambda  \int\limits_0^{r_0}e^{u_{\m}(r)}dr&=&\m \e_\lambda  \int\limits^{r_0}_{r_0-\delta}e^{u_{1}(r)}dr+\m \e_\lambda  \int\limits_{r_0-2\delta}^ {r_0- \delta}e^{u_{2}(r)}dr+\m \e_\lambda  \int\limits_0^{r_0-2\delta}e^{u_{3}(r)}dr \\
&= &\m \e^2_\lambda  \int\limits^{0}_{ -\delta/\eps}e^{u_{1}(\e_\lambda s+r_0)}dr+O\(\eps^{1+\sigma}\)\to  \int\limits_ {\rr}e^{w(s)}d s \quad  \hbox{  as}\ \lambda\rightarrow 0,
\end{eqnarray*}
because of \eqref{cc1}, \eqref{cc2} and dominate convergence Lebesgue's Theorem.
That proves \eqref{conv2}.

\end{Proof}

\section{A contraction mapping argument and the proof of the main theorem}\label{sec4}
First of all we point out  that $u_{\m}+\phi_{\m}$ is a solution to \eqref{P} if and only if
  $\phi_{\m}$ is a solution of the problem
\begin{equation}\label{L2}
\left\{
\begin{array}{lr}
\mathcal{L}_{\m}(\phi_{\m})=\mathcal{N}_{\m}(\phi_{\m})+\mathcal{R}_{\m}(u_{\m})\qquad \mbox{in}\,\, (0, r_0)\\\\
\phi_{\m}'(0)=\phi'_{\m}(r_0)=0
\end{array}
\right.
\end{equation}
where $R_\m(u_{\m})$ is given in \eqref{rlam},
 $$\mathcal{L}_{\m}(\phi_{\m}):=-\phi_{\m}''-\frac{N-1}{r}\phi_{\m}'+\phi_{\m}-\m e^{u_{\m}}\phi_{\m}$$
and
$$ \mathcal{N}_{\m}(\phi_{\m}):=\m e^{u_{\m}+\phi_{\m}}-\m e^{u_{\m}}-\m e^{u_{\m}}\phi_{\m}.$$

\medskip

The next result  state that the linearized operator $\mathcal{L}_{\m}$ is uniformly invertible.
\\
\begin{proposition}\label{inv}
There exists $\m_0>0$ and $C>0$ such that for any $\m \in(0, \m_0)$ and for any $h\in L^{\infty}((0, r_0))$ there exists a unique $\phi\in W^{2, 2}((0, r_0))$ solution of
\begin{equation}
\left\{
\begin{array}{lr}
\mathcal{L}_{\m}(\phi)=h\\
\phi'(0)=\phi'(r_0)=0
\end{array}
\right.
\end{equation}
which satisfies $$\|\phi\|_{L^\infty(0,r_0)}\leq C\|h\|_{L^1(0,r_0)}$$
\end{proposition}
\begin{Proof}
By contradiction we assume that there exist sequences $\m_n\rightarrow 0$, $h_n\in L^{\infty}((0, r_0))$ and $\phi_n\in W^{2, 2}((0, r_0))$ solutions of
\begin{equation}\label{lin}
\left\{
\begin{array}{lr}
-\phi_n''-\frac{N-1}{r}\phi_n'+\phi_n-\m_n e^{u_{\m_n}}\phi_n=h_n \qquad \mbox{in}\,\, (0, r_0)\\
\phi_n'(0)=\phi_n'(r_0)=0
\end{array}
\right.
\end{equation}
and
\begin{equation}\label{norme}
\|\phi_n\|_{L^\infty}=1\qquad \|h_n\|_{L^1}\rightarrow 0.
\end{equation}
Let $\psi_n(s)=\phi_n(\e_n s +r_0)$. Then $\psi_n$ solves
\begin{equation}\label{6.10}
\left\{
\begin{array}{lr}
-\psi_n''-\frac{N-1}{\e_n s +r_0}\e_n \psi_n'+\e_n^2 \psi_n -\m_n \e_n^2 e^{u_n(\e_n s+ r_0)} \psi_n=\e_n^2 h_n(\e_n s+r_0)\quad \mbox{in}\,\, \left(-\frac{r_0}{\e_n}, 0\right)\\
\psi_n'(-\frac{r_0}{\e_n})=\psi'_n(0)=0
\end{array}
\right.
\end{equation}
We point out that, since $\psi_n$ is bounded in $L^{\infty}((0, r_0))$, we get that, by standard elliptic regularity theory,
 $\psi_n \rightarrow \psi$ $C^2-$ uniformly on compact sets of $(-\infty, 0]$. \\
 Hence we multiply the equation in \eqref{6.10} by a $C^{\infty}_0$- test function, we integrate
 and we use  \eqref{conv1}
 to deduce that $\psi$ solves

\begin{equation}\label{6.11}
\left\{
\begin{array}{lr}
-\psi''-e^w \psi=0 \qquad \mbox{in}\,\, (-\infty, 0)\\
\|\psi\|_{\infty}\leq 1\\
\psi'(0)=0.
\end{array}
\right.
\end{equation}
A straightforward computation shows (see Lemma 4.2, \cite{g2}) that there exist $a,b\in\rr$ such that
 $$\psi(s)=a\frac{e^{\sqrt{2}s}-1}{e^{\sqrt{2}s}+1}+b\left(-2+\sqrt{2}s\frac{e^{\sqrt{2}s}-1}{e^{\sqrt{2}s}+1}\right).$$
 It is immediate to check that $b=0,$ since $\|\psi\|_{\infty}\leq 1$ and also that $a=0,$ since
   $\psi'(0)=0.$ Therefore, $\psi\equiv 0$ in $(0,r_0)$.\\
We claim that $\|\phi_n\|_{\infty}=o(1)$. This immediately gives a contradiction since by assumption $\|\phi_n\|_{\infty}=1$.
To prove the claim we introduce the function $G$ being
the Green function of the operator $-u''-\frac{N-1}{r}u'+u$ with Neumann boundary condition.\\
By \eqref{lin}, we deduce that
\begin{eqnarray*}
\phi_n(r)&=& \int_0^{r_0}G(r, t)\m_n e^{u_{\m_n}}\phi_n(t)\, dt +\int_0^{r_0}G(r, t)h_n(t)\, dt\\
&=&\e_n\m_n \int_{-\frac{r_0}{\e_n}}^0 G(r, \e_n s +r_0)e^{u_{\m_n}(\e_n s+r_0)}\psi_n(s)\, ds +\int_0^{r_0}G(r, t)h_n(t)\, dt\\
&=& G(r) \e_n\m_n \int_{-\frac{r_0}{\e_n}}^0 e^{u_{\m_n}(\e_n s+r_0)}\psi_n(s)\, ds +\int_0^{r_0}G(r, t)h_n(t)\, dt\\
&&+\e_n\m_n\int_{-\frac{r_0}{\e_n}}^0\left[G(r, \e_n s+ r_0)-G(r)\right]e^{u_{\m_n}(\e_n s+r_0)}\psi_n(s)\, ds
\end{eqnarray*}
Since $G$ is bounded, it is immediate to check that $\int_0^{r_0}G(r, t)h_n(t)\, dt=o(1)$. We want to show that also
\begin{equation}\label{claim1}
\e_n\m_n\int_{-\frac{r_0}{\e_n}}^0 \left[G(r, \e_n s+r_0)-G(r)\right]e^{u_{\m_n}(\e_ns+r_0)}\psi_n(s)\, ds=o(1)
\end{equation}
If this is true then $$\phi_n(r)=G(r) K_n +o(1)$$ where $$K_n:=\e_n\m_n \int^0_{-\frac{r_0}{\e_n}}e^{u_{\m_n}(\e_ns+r_0)}\psi_n(s)\, ds$$
We compute $$G(r_0)K_n +o(1)=\phi_n(r_0)=\psi_n(0)=o(1)$$ and hence $K_n=o(1)$ since $G(r_0)\neq 0$. Then
$\|\phi_n\|_{\infty}=o(1)$ and this gives a contradiction.\\
It remains to prove  \eqref{claim1}. We have:\\
\begin{eqnarray*}
&&\left|\e_n\m_n \int_{-\frac{r_0}{\e_n}}\left[G(r, \e_n s+r_0)-G(r)\right]e^{u_{\m_n}(\e_n s+r_0)}\psi_n(s)\, ds
\right|\leq \e_n^2\m_n\int_{-\frac{r_0}{\e_n}}^0
|s|e^{u_{\m_n}(\e_ns+r_0)}|\psi_n(s)|\, ds
\\
&&= \underbrace{\e_n^2\m_n\int_{-\frac{\delta_n}{\e_n}}^{0}|s|e^{u_{1_n}(\e_ns+r_0)}|\psi_n(s)|\, ds}_{(I)}+
+\underbrace{\e_n^2\m_n\int_{-\frac{2\delta_n}{\e_n}}^{-\frac{\delta_n}{\e_n}}
|s|e^{ u_{2_n}(\e_ns+r_0)}|\psi_n(s)|\, ds}_{(II)}\\
&&+\underbrace{\e_n^2\m_n\int_{-\frac{r_0}{\e_n}}^{-\frac{2\delta_n}{\e_n}}
|s|e^{u_{3_n}(\e_ns+r_0)}|\psi_n(s)|\, ds}_{(III)}=o(1)
\end{eqnarray*}
Indeed, taking into account  that $\psi_n\rightarrow 0$ pointwise in $(-\infty, 0)$ and $\|\psi_n\|_{\infty}\leq 1,$
 by \eqref{cc2} we deduce
\begin{equation*}
(I)=O\left(\int_{-\infty}^0|s| e^{-a |s|}|\psi_n(s)|\, ds\right)=o(1)
\end{equation*}
for some $a>0$, and arguing as in Step 2 and in Step 3 of Lemma \ref{lemerr}, we get  respectively
\begin{equation*}
(III)=O\left(\int_{-\frac{r_0}{\e_n}}^{-\frac{2\delta_n}{\e_n}}|s|e^{-|s|}|\psi_n(s)|\, ds\right)
=O\left(\int_{-\infty}^0 |s|e^{-|s|}|\psi_n(s)|\, ds\right)=o(1)
\end{equation*}
 \begin{equation*}
(II)=O\left(\int_{-\frac{2\delta_n}{\e_n}}^{-\frac{\delta_n}{\e_n}}|s|e^{-|s|}|\psi_n(s)|\, ds\right)
=O\left(\int_{-\infty}^0|s| e^{-|s|}|\psi_n(s)|\, ds\right)=o(1).
\end{equation*}

\end{Proof}

\medskip

Finally, we are in position to use a contraction mapping argument to prove Theorem \ref{principale}.  \\

\begin{Proof}[Proof of Theorem \ref{principale}]
By Proposition \ref{inv}, we deduce that the linear operator $\mathcal{L}_{\m}$ is uniformly invertible and so problem \eqref{L2} can be rewritten as
\begin{equation}\label{conmap}
\phi =\mathcal{T}_\m(\phi ):=\mathcal{L}_{\m}^{-1}\left[\mathcal{R}_{\m}(\bar u_\m)+\mathcal{N}_\m (\phi )\right].\end{equation}
For a given number $\rho>0$ let us
consider the closed set  $A_{\rho}:=\left\{\phi\in L^{\infty}((0,r_0))\ :\  \|\phi\|_{\infty}\leq \rho \e_\m^{1+\sigma}\right\}$  where
$\eps_\m$ is defined in \eqref{sceltae} and $\sigma>0$ is given in Lemma \ref{lemerr}.\\

We will prove that if $\lambda$ is small enough, then  $ \mathcal{T}_\m:A_{\rho}\to A_{\rho}$ is a contraction map.

First of all, by \eqref{conv2} we get
 $$\|\mathcal{N}_\m(\phi)\|_{L^1}
\leq \|\m e^{u_{\m}}\|_{L^1}\|\phi\|_{L^\infty}^2
\leq \frac{C}{\e_\m}\|\phi \|^2_{L^\infty}\quad \hbox{for any}\ \phi\in A_\rho$$
and also
$$\|\mathcal{N}_{\m}(\phi_1)-\mathcal{N}_{\m}(\phi_2)\|_{L^1}\leq \frac{C}{\e_\m}\left(\max_{i=1, 2}\|\phi_i\|_{L^\infty}\right)\|\phi_1-\phi_2\|_{L^\infty}\quad \hbox{for any}\ \phi_1,\phi_2\in A_\rho$$
for some $C>0$.\\
By Lemma \ref{lemerr} we deduce that for some $\rho>0$
$$\|\mathcal{T}_\m (\phi)\|_{L^\infty}\leq C\left(\|\mathcal{R}_{\m}(u_{\m})\|_{L^1}+\|\mathcal{N}_{\m}(\phi)\|_{L^1}\right)\leq \rho\e_\m^{1+\sigma}$$
 and so $ \mathcal{T}_\m$ maps $A_{\rho}$ into itself. Moreover
 $$\|\mathcal{T}_\m(\phi_1)-\mathcal{T}_\m(\phi_2)\|_{L^\infty}\leq C\|\mathcal{N}_{\m}(\phi_1)-\mathcal{N}_{\m}(\phi_2)\|_{L^1}\leq C\e_\lambda^{\sigma}\|\phi_1-\phi_2\|_{L^\infty}$$
 which proves that  for $\m$ small enough $\mathcal{T}_\m$ is a contraction mapping on $A_\rho,$ for a suitable $\rho$.\\

 Therefore, $\mathcal{T}_\m$ has a unique fixed point in  $A_\rho$, namely there exists a unique solution $\phi=\phi_\m\in A_\rho$ of the equation \eqref{conmap} or equivalently there exists a unique solution $u_\m+\phi_\m$ of problem \eqref{P}. \\
 Estimate \eqref{convergenza} follows by the definition of $u_\m$ which coincides with $u_3$ far away from $r_0.$ Indeed if $[a,b]$ is a compact set in
 $(0, r_0),$ we get that for $\m$ small enough
 $$\e _\m\(u_{\m}(r)+\phi_\m(r)\)= (A_1+A_2\e _\m+A_3\e _\m^2)\mathcal{U}(r)+\e _\m\phi_\m(r)\rightarrow \frac{\sqrt{2}}{\mathcal{U}'(r_0)}\mathcal{U}(r )\ \hbox{as}\ \lambda\to0,$$
because of  \eqref{abc} and the fact $ \|\phi_\m\|_{L^\infty}\to0$ as $\m\to0.$ \\
 Finally, estimate \eqref{convergenzaL1}  follows by \eqref{conv2}, taking into account that  $ \|\phi_\m\|_{L^\infty}\to0$ as $\m\to0.$

  \end{Proof}

\end{document}